\documentclass[a4paper,oneside,11pt]{article}

\usepackage[francais]{babel} 
\usepackage{vmargin}
\usepackage{epsfig}
\psfigdriver{dvips}
\usepackage{graphicx}
\usepackage{color}
\usepackage[latin1]{inputenc}
\usepackage{amsmath}
\usepackage{amssymb}
\usepackage[all]{xy}

\begin{document}
\renewcommand{\thepage}{\empty} 

\newtheorem{thm}{Théorème}
\newtheorem{cor}[thm]{Corollaire}
\newtheorem{lem}[thm]{Lemme}
\newtheorem{prop}[thm]{Proposition}
\newtheorem{dfn}[thm]{Définition}
\newtheorem{rem}[thm]{Remarque}
\newtheorem{exm}[thm]{Exemple}
\newcommand{\al}{\alpha}
\newcommand{\A}{\mathcal{A}}
\newcommand{\Aa}{\mathcal{A}^\alpha}
\newcommand{\ta}{\theta}

\sloppy 

\begin{center}
{\Large\bf Stabilit\'e de la propri\'et\'e de Koszul
 pour les alg\`ebres homog\`enes vis-\`a-vis du produit semi-crois\'e}
\end{center} 
\vspace{0,5cm}
\begin{center}
Antonin POTTIER
\end{center}
\begin{center}
\small {\'Ecole Normale Sup\'erieure, 45, rue d'Ulm, 75230 Paris Cedex 05 et \\Laboratoire de Physique Th\'eorique, UMR 8627, Universit\'e Paris XI,\\
B\^atiment 210, F-91 405 Orsay Cedex, France, 
antonin.pottier@ens.fr}
\end{center}

\begin{abstract}
We study the stability of Koszul and Gorentein properties for the semi-cross product of homogeneous algebras.\\
-----------------------------------
Nous étudions la conservation des propriétés de Koszul et de
Gorenstein pour le produit semi-croisé des algèbres homogènes.
\end{abstract}

\section{Introduction}

Le but de cette note est d'étudier la stabilité de certaines
propriétés homologiques des algèbres homogènes par produit
semi-croisé, introduit au paragraphe 7.1 de \cite{Moduli}. Plus
précisèment, nous montrons qu'une algèbre homogène est de type
Koszul si et seulement si un de ces produits semi-croisés l'est.
Dans le cas où la dimension globale est finie, être de type
Gorenstein est équivalent pour l'algèbre et ses produits
semi-croisés.

Différentes notions relatives aux algèbres quadratiques introduites par \cite{Priddy} sont
généralisées aux algèbres homogènes dans \cite{Homo}. En particulier un $N$--complexe est
canoniquement attaché à toute algèbre $N$--homogène, dont le complexe de Koszul de \cite{Berg} est
une contraction.
Dans
l'article \cite{Berg}, il est montré qu'être de type Koszul pour une
algèbre homogène est équivalent à l'acyclicité de ce
complexe. C'est cette caractérisation que nous utiliserons. En plus
d'algèbres quadratiques, on trouve des algèbres cubiques dans la classification
des algèbres régulières de dimension $3$ décrite par \cite{AS}. D'autres exemples
d'algèbres homogènes de degré supérieur à $3$ ont été étudiées par la suite dans \cite{Berg}, ainsi
que dans \cite{YM} et \cite{YM2} en liaison avec certaines équations issues de la physique théorique.


\section{Rappels et notations}

$k$ est un corps fixé dans toute la suite, tous les produits
tensoriels seront pris sur $k$, $\otimes=\otimes_k$. Soit
$\A=A(E,R)$ une algèbre homogène de degré $N$. C'est le quotient de
l'algèbre tensorielle $T(E)$ associée à un $k$--espace vectoriel $E$
de dimension finie par un idéal bilatère $I(R)$ engendré par un
espace de relations $R\subset E^{\otimes N}$. Soit $\alpha$ un
automorphisme de l'algèbre graduée $\A$. Il est défini par un
automorphisme de $E$ étendu canoniquement à $T(E)$, encore noté
$\alpha$ et tel que $\alpha(R)=R$. L'algèbre $\Aa$, produit
semi-croisé de $\A$ par $\alpha$, est donnée par l'espace vectoriel
gradué sous-jacent à $\A$ muni du produit $\cdot$ défini sur les
éléments homogènes par $x\cdot y=x\al^{|x|}(y)$ où $|x|$ est le
degré de $x$ et où le symbole pour le produit dans $\A$ est omis,
voir \cite{Moduli}. $\Aa$ est encore une algèbre associative avec
unité, identique à celle de $\A$. Remarquons tout de suite que
$id:\Aa\rightarrow\A$ est un isomorphisme de $k$--espaces vectoriels
et que $\alpha$ est encore un automorphisme de l'algèbre $\Aa$,
ainsi que $\A=(\Aa)^{\al^{-1}}$.

Définissons maintenant $\ta$, automorphisme de l'espace vectoriel
gradué $T(E)$: en degré $n+1$, $\ta_{n+1}(x_0\otimes x_1\otimes
\ldots\otimes x_n)=x_0\otimes \al(x_1)\otimes \ldots\otimes
\al^n(x_n)$. Relativement à la décomposition $E^{\otimes n+1}\simeq
E^{\otimes p+1}\otimes E^{\otimes n-p}$, on a la formule:
\begin{equation} \ta_{n+1}=(\ta_{p+1}\otimes id)
\circ(id\otimes(\alpha^{p+1}\circ\ta_{n-p}))
\label{dec}\end{equation}

Comme application de ces définitions, prouvons la proposition
suivante.

\begin{prop} $\Aa$ est une algèbre homogène de degré $N$,
$\Aa=A(E,\ta_N^{-1}(R))$.

\end{prop}

Considérons $m:T(E)\rightarrow\A$ défini en degré $n$ par
$m(x_1\otimes x_2\otimes \ldots\otimes x_n)=x_1x_2\ldots x_n$ et
$m_\al:T(E)\rightarrow\Aa$ défini en degré $n$ par $m_\al(x_1\otimes
x_2\otimes \ldots\otimes x_n)=x_1\cdot x_2\cdot\ldots \cdot x_n$.
Alors le diagramme suivant d'applications linéaires commute,
c'est-à-dire $m_\al=m\circ \ta$.

$$\xymatrix{\relax
     T(E) \ar[r]^-{m_\al} \ar[d]^\ta & \Aa \ar[d]^{id} \\
     T(E) \ar[r]^-{m}  & \A
}$$

Par définition \cite{Homo}, $Ker\, m = I(R)$, d'où l'égalité $Ker\,
m_\al=\ta^{-1}(I(R))=I(\ta_N^{-1}(R))$. En conséquence $\Aa$ est une
algèbre homogène de degré $N$, $\Aa=A(E,\ta_N^{-1}(R))$.

\medskip

\textbf{Exemple} Soient $\A=A(E=kx\oplus ky,x\otimes y\otimes
x-y\otimes x\otimes y)$ l'algèbre des tresses à $3$ brins et $\al$
l'automorphisme involutif échangeant $x$ et $y$. Alors le produit
semi-croisé de $\A$ par $\al$ est $\Aa=A(E,x\otimes x\otimes
x-y\otimes y\otimes y)$, ce qui est une écriture plus symétrique. Nous poursuivrons
plus loin l'étude de cette algèbre via son produit semi-croisé.

\section{Conservation des types Koszul et Gorenstein}
\begin{thm} $\A$ est de type Koszul si et seulement si $\Aa$ est de
type Koszul.
\end{thm}

D'après \cite{Homo}, $\A$ est de type Koszul si le complexe
$\mathcal{C}$ de $\A$--modules à gauche est acyclique en degrés strictement positifs. Le complexe
$\mathcal{C}$, c'est-à-dire \xymatrix{\relax
     \ldots \ar[r]^-{d} & \A\otimes\A^{!*}_{(p+1)N} \ar[r]^-{\delta}
     & \A\otimes\A^{!*}_{pN+1} \ar[r]^-{d}&\A\otimes\A^{!*}_{pN} \ar[r]^-{\delta}&
   \ldots  \ar[r]^-{d}&\A \ar[r] &0  \\
}, est la contraction $\mathcal{C}=C_{N-1,0}$ du $N$--complexe
$K(\A)$ de $\A$--modules à gauche (avec $\delta=d^{N-1}$) \xymatrix{
     \ldots \ar[r]^-{d} & \A\otimes\A^{!*}_{i+1} \ar[r]^-{d}
     & \A\otimes\A^{!*}_{i} \ar[r]^-{d}&\A\otimes\A^{!*}_{i-1} \ar[r]^-{d}&
   \ldots  \ar[r]^-{d}&\A \ar[r] &0  \\
}.

 Nous pouvons voir $(K(\A),d)$ comme un $N$--complexe d'espaces
vectoriels. Nous allons construire un isomorphisme de $N$--complexes
entre $(K(\Aa),d^\al)$ et $(K(\A),d)$. Cela induira un isomorphisme
de complexes entre leur contraction. Un isomorphisme de complexes
étant un homologisme, l'acyclicité de $\mathcal{C}^\al$ sera
équivalente à celle de $\mathcal{C}$, ce qui prouvera le théorème.

Rappelons que $\A^{!*}_i$ est naturellement un sous-espace de
$E^{\otimes i}$ (cf. \cite{Homo}). Définissons
$K(\ta):K(\Aa)\rightarrow K(\A)$ en degré $i$ par:

$$K(\ta)_i: K(\Aa)_i=\Aa\otimes (\Aa)  _{i} ^{!*} \rightarrow
K(\A)_i=\A\otimes\A^{!*}_{i} $$\begin{equation} a \otimes e \mapsto
a \otimes \al^{|a|}\circ \ta_i(e)\end{equation}

Il est clair que $K(\ta)_i$ est un isomorphisme d'espaces
vectoriels. Vérifions alors que $K(\ta)$ est un morphisme de
$N$--complexes.
$$
\xymatrix{\relax
     \ldots \ar[r]^-{d^\al} & \Aa\otimes(\Aa)^{!*}_{i+1}
     \ar[r]^-{d^\al} \ar[d]^{K(\ta)_{i+1}}
     & \Aa\otimes(\Aa)^{!*}_{i} \ar[r]^-{d^\al} \ar[d]^{K(\ta)_{i}}
     &\Aa\otimes(\Aa)^{!*}_{i-1} \ar[r]^-{d^\al} \ar[d]^{K(\ta)_{i-1}}&
   \ldots   \\
\ldots \ar[r]^-{d} & \A\otimes\A^{!*}_{i+1} \ar[r]^-{d}
     & \A\otimes\A^{!*}_{i} \ar[r]^-{d}&\A\otimes\A^{!*}_{i-1} \ar[r]^-{d}&
   \ldots   \\
}$$

 Soit $a\otimes (e\otimes f)$ un élément générique de
$\Aa\otimes(\Aa)_{i+1}^{!*}$ avec $\A^{!*}_{i+1}\subset E^{\otimes
i+1}\simeq E\otimes E^{\otimes i}$. D'une part $d^\al(a\otimes
(e\otimes f))=a\cdot e \otimes f=a \al^{|a|}(e)\otimes f$, donc
$K(\ta)_i\circ d^\al(a\otimes (e\otimes f))=a \al^{|a|}(e)\otimes
\al^{|a|+1}(\ta_i(f))$ car $|a \al^{|a|}(e)|=|a|+1$ puisque $\al$
est de degré $0$. D'autre part $K(\ta)_{i+1}(a\otimes (e\otimes
f))=a\otimes\al^{|a|}\circ \ta_{i+1}(e\otimes f)=a\otimes\al^{|a|}(
e\otimes \al\circ\ta_{i}(
f))=a\otimes(\al^{|a|}(e)\otimes\al^{|a|+1}(\ta_i(f)) )$ en
utilisant (\ref{dec}), donc $d\circ K(\ta)_{i+1}(a\otimes (e\otimes
f))=a\al^{|a|}(e)\otimes\al^{|a|+1}(\ta_i(f)) $. Finalement
$K(\ta)_i\circ d^\al=d\circ K(\ta)_{i+1}$, et $K(\ta)$ est un
isomorphisme de $N$--complexes d'espaces vectoriels. CQFD.

\begin{prop} Si $\A$ est de type Koszul de dimension globale finie alors $\Aa$
l'est aussi. \end{prop} En effet, dans le cas où $\A$ est de type
Koszul la dimension globale $D$ est donnée par le plus grand entier
tel que $\mathcal{C}_D\neq 0$ (avec $\mathcal{C}=\A$). Via
l'isomorphisme $K(\ta)$, $\mathcal{C}_D\neq 0$ équivaut à
$\mathcal{C}^\al_D\neq 0$, d'où la proposition.

\medskip

\textbf{Exemple} Montrons que $\Aa=A(E,R=x\otimes x\otimes
x-y\otimes y\otimes y)$ est de type Koszul de dimension globale $2$,
 ce qui montrera en vertu des théorèmes précédents que l'algèbre
des tresses à $3$ brins est du même type.

Le $3$--complexe $K(\Aa)$ se calcule simplement:
$$
\xymatrix{\relax
     0 \ar[r]
     & \Aa \ar[r]^-{d_3}
     & (\Aa)^4 \ar[r]^-{d_2}
     & (\Aa)^2 \ar[r]^-{d_1}
     & \Aa \ar[r]
     & 0 \\
}
$$
avec $d_3: a\mapsto (ax,0,0,-ay)$, $d_2:(a,b,c,d)\mapsto (ax+cy,bx+dy)$ et
$d_1: (a,b)\mapsto ax+by $.
Le complexe de Koszul $C^\al$ obtenu en contractant s'écrit, dans ce cas:
$$
\xymatrix{\relax
     0 \ar[r]& \Aa \ar[r]^-{\delta} & (\Aa)^2
     \ar[r]^-{d}
     & \Aa \ar[r]
     & 0    \\
}$$

 où $\delta^\al=d_2\circ d_3: a \mapsto a(x^2,-y^2)$ et $d=d_1$. La suite est exacte au niveau de $(\Aa)^2$,
c'est la définition de l'algèbre par générateurs et relations \cite{Homo}. Il suffit donc de vérifier
l'injectivité de la première flèche.

\emph{\underline{Lemme}:} $x$ et $y$ sont  réguliers à droite
 (sans diviseur de zéro à gauche).

 Raisonnons
par récurrence sur le degré des éléments de l'algèbre, l'initialisation étant évidente. Supposons $x$ et $y$ réguliers
jusqu'au degré $n$. Soit $a\in \Aa_{n+1}$ tel que $ax=0$, alors $d(a,0)=0$.
 De par l'exactitude au niveau de $(\Aa)^2$, il existe $a'\in \Aa_{n-1}$ tel que $\delta(a')=(a,0)$. Donc $a'y^2=0$ et
 par hypothèse de récurrence $a'=0$ d'où $a=a'x^2=0$.
  $x$ est bien régulier à droite jusqu'au degré $n+1$,
 la démonstration pour $y$ est identique. Le lemme est prouvé.

Puisque $x$ et $y$ sont réguliers à droite, la première flèche du complexe de Koszul $C^\al$ est
donc injective. Donc l'algèbre des tresses à $3$ brins possède la propriété de Koszul
et est de dimension globale $2$.

\emph{\underline{Remarque}:} La propriété de Koszul permet de calculer la série de
 Poincaré $P_\A(t)=\sum \mathrm{dim}(\A_n)t^n$ de $\A$.
 En effet d'après \cite{Popov}, on a la relation suivante:

 \begin{equation}
P_\A(t)\left(\sum_n \mathrm{dim}(\A^!_{Nn})t^{Nn}-\mathrm{dim}(\A^!_{Nn+1})t^{Nn+1}\right)=1
 \end{equation}

 Dans notre cas $N=3$, cela donne $1/P_\A(t)=1-2t+t^3=(1-t)(1-t-t^2)$. Ainsi l'algèbre des tresses à $3$ brins
 est à croissance exponentielle.

\begin{thm} Si $\A$ est de type Koszul de dimension
globale  finie $D$, alors $\A$ est de type Gorenstein si et
seulement si $\Aa$ l'est.\end{thm}

Dans les hypothèses du théorème, $\A$ est de type Gorenstein si la
cohomologie du complexe dual $\mathcal{C}'$ est nulle en degré
strictement inférieur à $D$. Le complexe de cochaînes $\mathcal{C}'$
de $\A$--modules à droite est obtenu à partir du complexe de chaînes
$\mathcal{C}$ de $\A$--modules à gauche en appliquant le foncteur
contravariant $Hom_\A(\bullet,\A)$. Le complexe de cochaînes
$\mathcal{C}'$ est la contraction $C_{1,0}$ du $N$--complexe $L(\A)$
obtenu en appliquant le foncteur contravariant $Hom_\A(\bullet,\A)$
à $K(\A)$ comme expliqué dans \cite{YM}. Or il est immédiat que
$Hom_\A(K(\ta),\A)$ est toujours un isomorphisme de $N$--complexes
d'espaces vectoriels entre $L(\Aa)$ et $L(\A)$. Il induit donc un
isomorphisme de complexes entre $(\mathcal{C}^\al)' $ et
$\mathcal{C}'$, d'où un homologisme, ce qui prouve le théorème.

\textbf{(Contre)-exemple} Dans le cas de l'algèbre des tresses à $3$ brins,
 le $3$--complexe $L(\Aa)$ s'obtient facilement à partir de $K(\Aa)$:

$$
\xymatrix{\relax
     0 \ar[r]& \Aa \ar[r]^-{d^1}
     & (\Aa)^2 \ar[r]^-{d^2}
     & (\Aa)^4 \ar[r]^-{d^3}
     & \Aa \ar[r]
     & 0    \\
}$$

avec $d^1: a\mapsto (xa,ya) $, $d^2:(a,b)\mapsto (xa,xb,ya,yb)$ et $d^3: (a,b,c,d)\mapsto xa-yd$.
 Le complexe de Gorenstein $(C^\al)'$ obtenu en contractant s'écrit donc, dans ce cas:
$$
\xymatrix{\relax
     0 \ar[r]& \Aa \ar[r]^-{d'}
& (\Aa)^2 \ar[r]^-{\delta '}
     & \Aa \ar[r]
     & 0    \\
}$$

où $d'=d^1$ et $\delta'=d^3\circ d^2: a \mapsto (x^2,-y^2)a$. Il est alors clair que $H^2((C^\al)')\neq k$,
 donc l'algèbre n'est pas de type Gorenstein.
  En résumé, l'algèbre des tresses à $3$--brins est Koszul de dimension $2$, mais n'est pas Gorenstein.


\end{document}